# Regression tree models for designed experiments*

**Wei-Yin Loh[1]**

*University of Wisconsin, Madison*

**Abstract:** Although regression trees were originally designed for large datasets, they can profitably be used on small datasets as well, including those from replicated or unreplicated complete factorial experiments. We show that in the latter situations, regression tree models can provide simpler and more intuitive interpretations of interaction effects as differences between conditional main effects. We present simulation results to verify that the models can yield lower prediction mean squared errors than the traditional techniques. The tree models span a wide range of sophistication, from piecewise constant to piecewise simple and multiple linear, and from least squares to Poisson and logistic regression.

## 1. Introduction

Experiments are often conducted to determine if changing the values of certain variables leads to worthwhile improvements in the mean yield of a process or system. Another common goal is estimation of the mean yield at given experimental conditions. In practice, both goals can be attained by fitting an accurate and interpretable model to the data. Accuracy may be measured, for example, in terms of prediction mean squared error, PMSE $= \sum_i E(\hat{\mu}_i - \mu_i)^2$, where $\mu_i$ and $\hat{\mu}_i$ denote the true mean yield and its estimated value, respectively, at the $i$th design point.

We will restrict our discussion here to complete factorial designs that are unreplicated or are equally replicated. For a replicated experiment, the standard analysis approach based on significance tests goes as follows. (i) Fit a full ANOVA model containing all main effects and interactions. (ii) Estimate the error variance $\sigma^2$ and use $t$-intervals to identify the statistically significant effects. (iii) Select as the "best" model the one containing only the significant effects.

There are two ways to control a given level of significance $\alpha$: the individual error rate (IER) and the experimentwise error rate (EER) (Wu and Hamda [22, p. 132]). Under IER, each $t$-interval is constructed to have individual confidence level $1 - \alpha$. As a result, if all the effects are null (i.e., their true values are zero), the probability of concluding at least one effect to be non-null tends to exceed $\alpha$. Under EER, this probability is at most $\alpha$. It is achieved by increasing the lengths of the $t$-intervals so that their simultaneous probability of a Type I error is bounded by $\alpha$. The appropriate interval lengths can be determined from the studentized maximum modulus distribution if an estimate of $\sigma$ is available. Because EER is more conservative than IER, the former has a higher probability of discovering the

*This material is based upon work partially supported by the National Science Foundation under grant DMS-0402470 and by the U.S. Army Research Laboratory and the U.S. Army Research Office under grant W911NF-05-1-0047.

[1]Department of Statistics, 1300 University Avenue, University of Wisconsin, Madison, WI 53706, USA, e-mail: loh@stat.wisc.edu

right model in the null situation where no variable has any effect on the yield. On the other hand, if there are one or more non-null effects, the IER method has a higher probability of finding them. To render the two methods more comparable in the examples to follow, we will use $\alpha = 0.05$ for IER and $\alpha = 0.1$ for EER.

Another standard approach is AIC, which selects the model that minimizes the criterion $\text{AIC} = n \log(\tilde{\sigma}^2) + 2\nu$. Here $\tilde{\sigma}$ is the maximum likelihood estimate of $\sigma$ for the model under consideration, $\nu$ is the number of estimated parameters, and $n$ is the number of observations. Unlike IER and EER, which focus on statistical significance, AIC aims to minimize PMSE. This is because $\tilde{\sigma}^2$ is an estimate of the residual mean squared error. The term $2\nu$ discourages over-fitting by penalizing model complexity. Although AIC can be used on any given collection of models, it is typically applied in a stepwise fashion to a set of *hierarchical* ANOVA models. Such models contain an interaction term only if all its lower-order effects are also included. We use the R implementation of stepwise AIC [14] in our examples, with initial model the one containing all the main effects.

We propose a new approach that uses a recursive partitioning algorithm to produce a set of nested piecewise linear models and then employs cross-validation to select a parsimonious one. For maximum interpretability, the linear model in each partition is constrained to contain main effect terms at most. Curvature and interaction effects are captured by the partitioning conditions. This forces interaction effects to be expressed and interpreted naturally—as contrasts of conditional main effects.

Our approach applies to unreplicated complete factorial experiments too. Quite often, two-level factorials are performed without replications to save time or to reduce cost. But because there is no unbiased estimate of $\sigma^2$, procedures that rely on statistical significance cannot be applied. Current practice typically invokes empirical principles such as *hierarchical ordering, effect sparsity*, and *effect heredity* [22, p. 112] to guide and limit model search. The hierarchical ordering principle states that high-order effects tend to be smaller in magnitude than low-order effects. This allows $\sigma^2$ to be estimated by pooling estimates of high-order interactions, but it leaves open the question of how many interactions to pool. The effect sparsity principle states that usually there are only a few significant effects [2]. Therefore the smaller estimated effects can be used to estimate $\sigma^2$. The difficulty is that a good guess of the actual number of significant effects is needed. Finally, the effect heredity principle is used to restrict the model search space to hierarchical models.

We will use the GUIDE [18] and LOTUS [5] algorithms to construct our piecewise linear models. Section 2 gives a brief overview of GUIDE in the context of earlier regression tree algorithms. Sections 3 and 4 illustrate its use in replicated and unreplicated two-level experiments, respectively, and present simulation results to demonstrate the effectiveness of the approach. Sections 5 and 6 extend it to Poisson and logistic regression problems, and Section 7 concludes with some suggestions for future research.

## 2. Overview of regression tree algorithms

GUIDE is an algorithm for constructing piecewise linear regression models. Each piece in such a model corresponds to a partition of the data and the sample space of the form $X \leq c$ (if $X$ is numerically ordered) or $X \in A$ (if $X$ is unordered). Partitioning is carried out recursively, beginning with the whole dataset, and the set of partitions is presented as a binary decision tree. The idea of recursive parti-tioning was first introduced in the AID algorithm [20]. It became popular after the appearance of CART [3] and C4.5 [21], the latter being for classification only.

CART contains several significant improvements over AID, but they both share some undesirable properties. First, the models are piecewise constant. As a result, they tend to have lower prediction accuracy than many other regression models, including ordinary multiple linear regression [3, p. 264]. In addition, the piecewise constant trees tend to be large and hence cumbersome to interpret. More importantly, AID and CART have an inherent bias in the variables they choose to form the partitions. Specifically, variables with more splits are more likely to be chosen than variables with fewer splits. This selection bias, intrinsic to all algorithms based on optimization through greedy search, effectively removes much of the advantage and appeal of a regression tree model, because it casts doubt upon inferences drawn from the tree structure. Finally, the greedy search approach is computationally impractical to extend beyond piecewise constant models, especially for large datasets.

GUIDE was designed to solve both the computational and the selection bias problems of AID and CART. It does this by breaking the task of finding a split into two steps: first find the variable $X$ and then find the split values $c$ or $A$ that most reduces the total residual sum of squares of the two subnodes. The computational savings from this strategy are clear, because the search for $c$ or $A$ is skipped for all except the selected $X$.

To solve the selection bias problem, GUIDE uses significance tests to assess the fit of each $X$ variable at each node of the tree. Specifically, the values (grouped if necessary) of each $X$ are cross-tabulated with the signs of the linear model residuals and a chi-squared contingency table test is performed. The variable with the smallest chi-squared $p$-value is chosen to split the node. This is based on the expectation that any effects of $X$ not captured by the fitted linear model would produce a small chi-squared $p$-value, and hence identify $X$ as a candidate for splitting. On the other hand, if $X$ is independent of the residuals, its chi-squared $p$-value would be approximately uniformly distributed on the unit interval.

If a constant model is fitted to the node and if all the $X$ variables are independent of the response, each will have the same chance of being selected. Thus there is no selection bias. On the other hand, if the model is linear in some predictors, the latter will have zero correlation with the residuals. This tends to inflate their chi-squared $p$-values and produce a bias in favor of the non-regressor variables. GUIDE solves this problem by using the bootstrap to shrink the $p$-values that are so inflated. It also performs additional chi-squared tests to detect local interactions between pairs of variables. After splitting stops, GUIDE employs CART's pruning technique to obtain a nested sequence of piecewise linear models and then chooses the tree with the smallest cross-validation estimate of PMSE. We refer the reader to Loh [18] for the details. Note that the use of residuals for split selection paves the way for extensions of the approach to piecewise nonlinear and non-Gaussian models, such as logistic [5], Poisson [6], and quantile [7] regression trees.

## 3. Replicated $2^4$ experiments

In this and the next section, we adopt the usual convention of letting capital letters $A$, $B$, $C$, etc., denote the names of variables as well as their main effects, and $AB$, $ABC$, etc., denote interaction effects. The levels of each factor are indicated in two ways, either by "$-$" and "$+$" signs, or as $-1$ and $+1$. In the latter notation, the variables $A$, $B$, $C$, $\ldots$, are denoted by $x_1$, $x_2$, $x_3$, $\ldots$, respectively.

TABLE 1
*Estimated coefficients and standard errors for $2^4$ experiment*

| | Estimate | Std. error | t | Pr(>\|t\|) |
|---|---|---|---|---|
| Intercept | 14.161250 | 0.049744 | 284.683 | < 2e-16 |
| x1 | -0.038729 | 0.049744 | -0.779 | 0.438529 |
| x2 | 0.086271 | 0.049744 | 1.734 | 0.086717 |
| x3 | -0.038708 | 0.049744 | -0.778 | 0.438774 |
| x4 | 0.245021 | 0.049744 | 4.926 | 4.45e-06 |
| x1:x2 | 0.003708 | 0.049744 | 0.075 | 0.940760 |
| x1:x3 | -0.046229 | 0.049744 | -0.929 | 0.355507 |
| x1:x4 | -0.025000 | 0.049744 | -0.503 | 0.616644 |
| x2:x3 | 0.028771 | 0.049744 | 0.578 | 0.564633 |
| x2:x4 | -0.015042 | 0.049744 | -0.302 | 0.763145 |
| x3:x4 | -0.172521 | 0.049744 | -3.468 | 0.000846 |
| x1:x2:x3 | 0.048750 | 0.049744 | 0.980 | 0.330031 |
| x1:x2:x4 | 0.012521 | 0.049744 | 0.252 | 0.801914 |
| x1:x3:x4 | -0.015000 | 0.049744 | -0.302 | 0.763782 |
| x2:x3:x4 | 0.054958 | 0.049744 | 1.105 | 0.272547 |
| x1:x2:x3:x4 | 0.009979 | 0.049744 | 0.201 | 0.841512 |

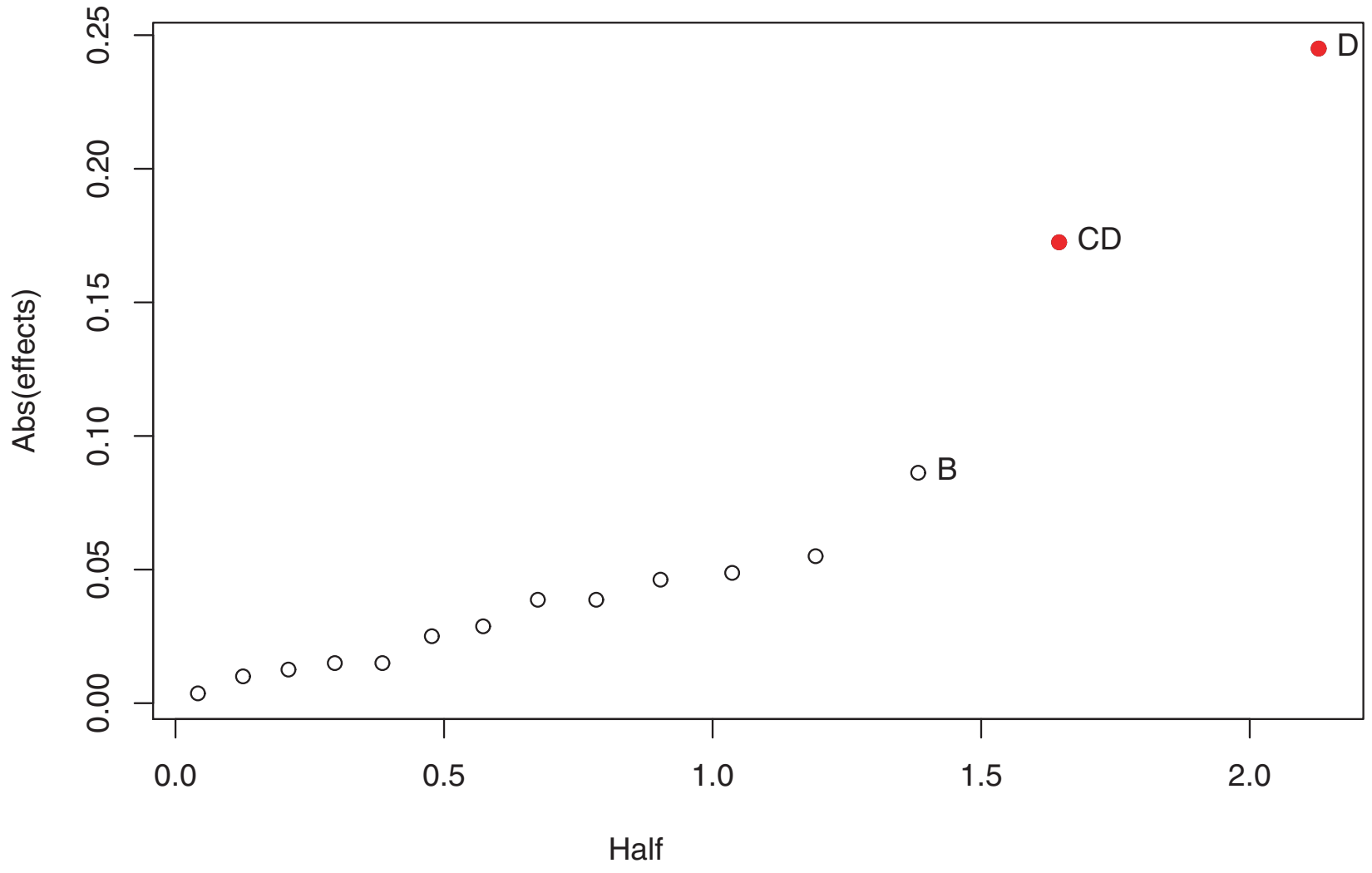

FIG 1. *Half-normal quantile plot of estimated effects from replicated $2^4$ silicon wafer experiment.*

We begin with an example from Wu and Hamada [22, p. 97] of a $2^4$ experiment on the growth of epitaxial layers on polished silicon wafers during the fabrication of integrated circuit devices. The experiment was replicated six times and a full model fitted to the data yields the results in Table 1.

Clearly, at the 0.05-level, the IER method finds only two statistically significant effects, namely $D$ and $CD$. This yields the model

$$\hat{y} = 14.16125 + 0.24502x_4 - 0.17252x_3x_4 \tag{3.1}$$

which coincides with that obtained by the EER method at level 0.1.

Figure 1 shows a half-normal quantile plot of the estimated effects. The $D$ and $CD$ effects clearly stand out from the rest. There is a hint of a $B$ main effect, but it is not included in model (3.1) because its $p$-value is not small enough. The $B$ effect appears, however, in the AIC model

$$\hat{y} = 14.16125 + 0.08627x_2 - 0.03871x_3 + 0.24502x_4 - 0.17252x_3x_4. \tag{3.2}$$

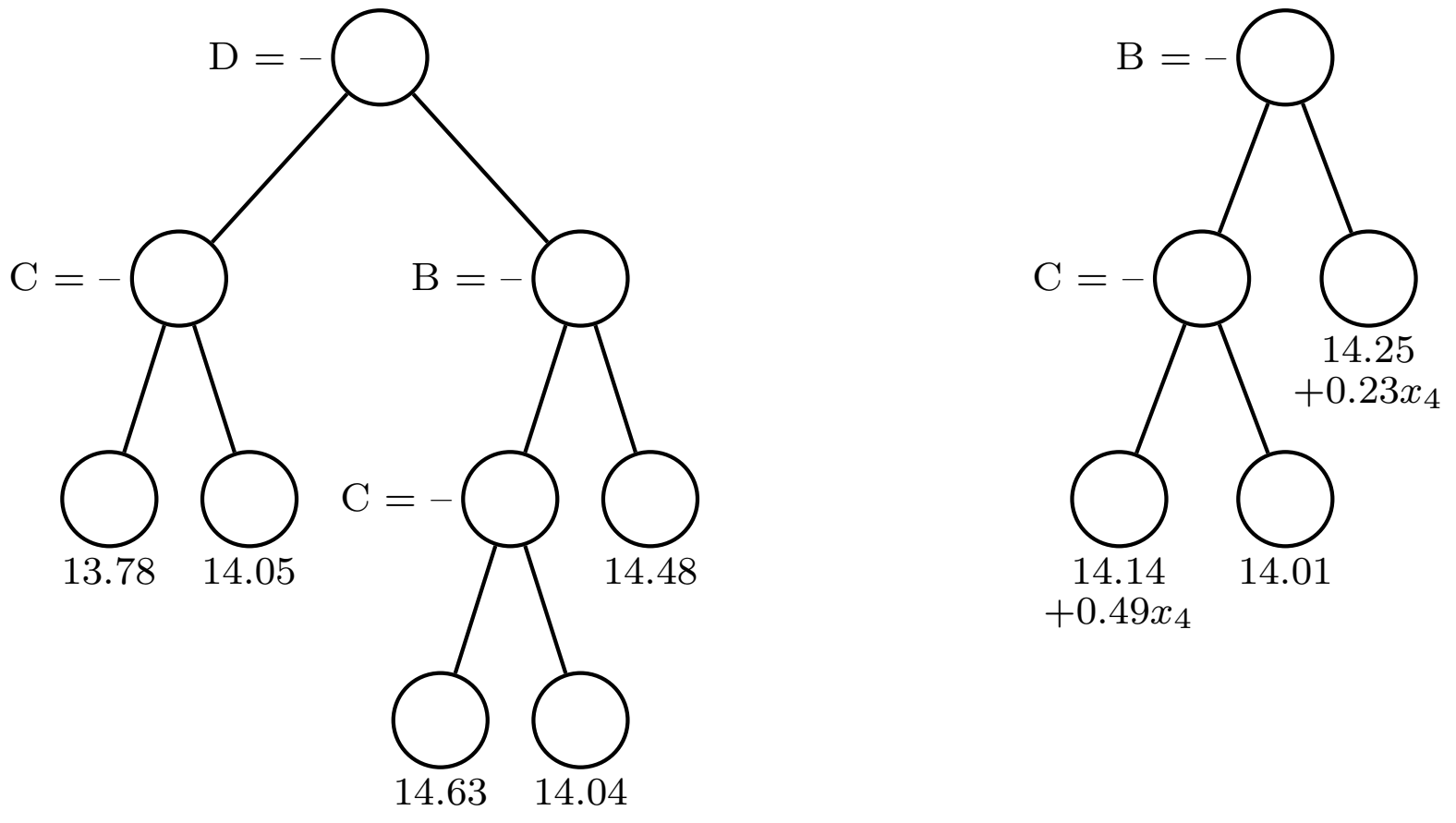

FIG 2. *Piecewise constant (left) and piecewise best simple linear or stepwise linear (right) GUIDE models for silicon wafer experiment. At each intermediate node, an observation goes to the left branch if the stated condition is satisfied; otherwise it goes to the right branch. The fitted model is printed beneath each leaf node.*

Note the presence of the small $C$ main effect. It is due to the presence of the $CD$ effect and to the requirement that the model be hierarchical.

The piecewise constant GUIDE tree is shown on the left side of Figure 2. It has five leaf nodes, splitting first on $D$, the variable with the largest main effect. If $D = +$, it splits further on $B$ and $C$. Otherwise, if $D = -$, it splits once on $C$. We observe from the node sample means that the highest predicted yield occurs when $B = C = -$ and $D = +$. This agrees with the prediction of model (3.1) but not (3.2), which prescribes the condition $B = D = +$ and $C = -$. The difference in the two predicted yields is very small though. For comparison with (3.1) and (3.2), note that the GUIDE model can be expressed algebraically as

$$
\begin{aligned}
\hat{y} &= 13.78242(1-x_4)(1-x_3)/4 + 14.05(1-x_4)(1+x_3)/4 \\
&\quad + 14.63(1+x_4)(1-x_2)(1-x_3)/8 + 14.4775(1+x_4)(1+x_2)/4 \\
&\quad + 14.0401(1+x_4)(1-x_2)(1+x_3)/8 \\
&= 14.16125 + 0.24502x_4 - 0.14064x_3x_4 - 0.00683x_3 \\
&\quad + 0.03561x_2(x_4+1) + 0.07374x_2x_3(x_4+1).
\end{aligned}
\tag{3.3}
$$

The piecewise best simple linear GUIDE tree is shown on the right side of Figure 2. Here, the data in each node are fitted with a simple linear regression model, using the $X$ variable that yields the smallest residual mean squared error, provided a statistically significant $X$ exists. If there is no significant $X$, i.e., none with absolute $t$-statistic greater than 2, a constant model is fitted to the data in the node. In this tree, factor $B$ is selected to split the root node because it has the smallest chi-squared $p$-value after allowing for the effect of the best linear predictor. Unlike the piecewise constant model, which uses the variable with the largest main effect to split a node, the piecewise linear model tries to keep that variable as a linear predictor. This explains why $D$ is the linear predictor in two of the three leaf nodes of the tree. The piecewise best simple linear GUIDE model can be expressed as

$$
\begin{aligned}
\hat{y} &= (14.14246 + 0.4875417x_4)(1 - x_2)(1 - x_3)/4 \\
&\quad + 14.0075(1 - x_2)(1 + x_3)/4 \\
&\quad + (14.24752 + 0.2299792x_4)(1 + x_2)/2 \\
&= 14.16125 + 0.23688x_4 + 0.12189x_3x_4(x_2 - 1) \\
&\quad + 0.08627x_2 + 0.03374x_3(x_2 - 1) - 0.00690x_2x_4.
\end{aligned} \tag{3.4}
$$

Figure 3, which superimposes the fitted functions from the three leaf nodes, offers a more vivid way to understand the interactions. It shows that changing the level of $D$ from $-$ to $+$ never decreases the predicted mean yield and that the latter varies less if $D = -$ than if $D = +$. The same tree model is obtained if we fit a piecewise multiple linear GUIDE model using forward and backward stepwise regression to select variables in each node.

A simulation experiment was carried out to compare the PMSE of the methods. Four models were employed, as shown in Table 2. Instead of performing the simula-

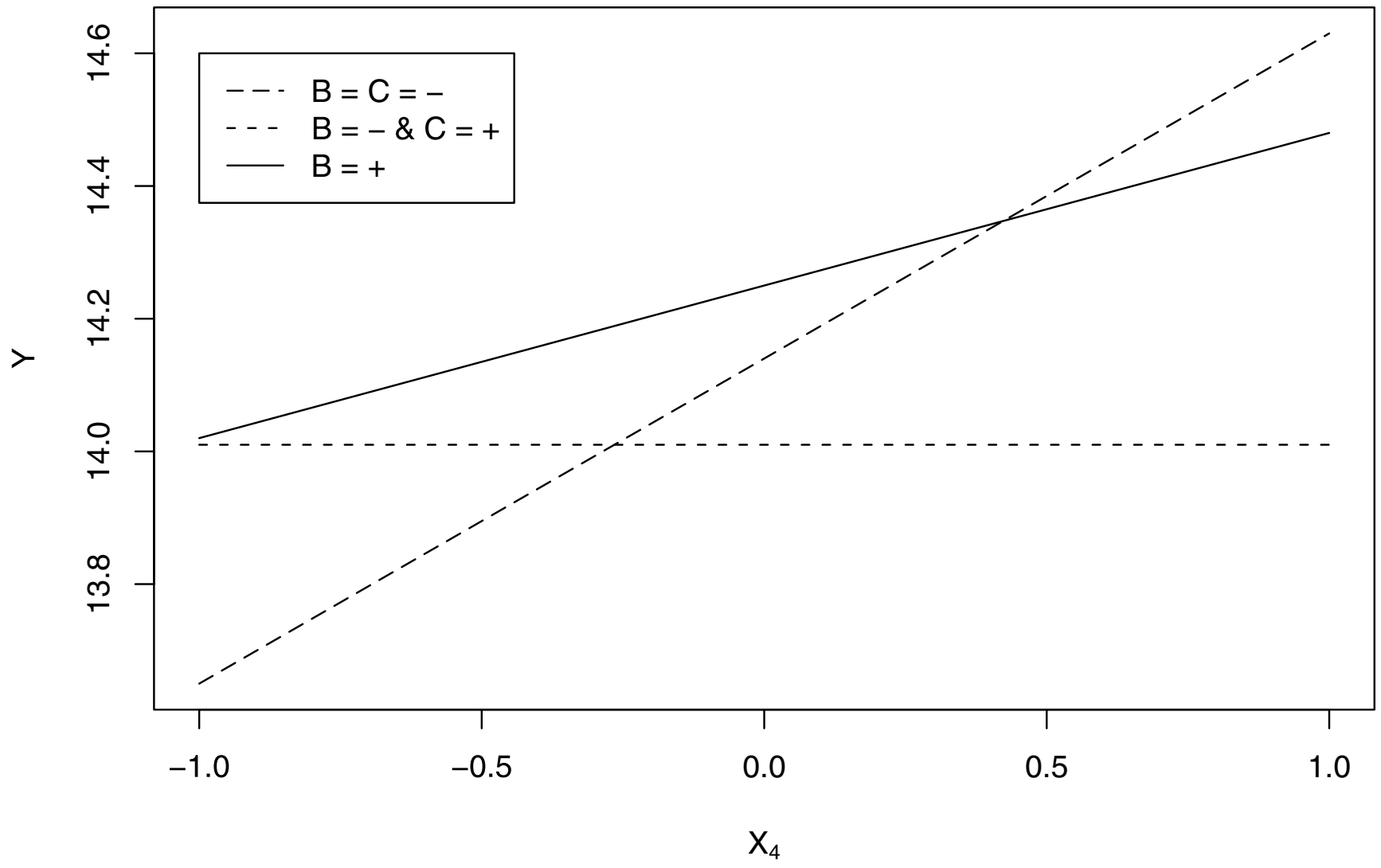

FIG 3. *Fitted values versus $x_4$ (D) for the piecewise simple linear GUIDE model shown on the right side of Figure 2.*

TABLE 2

*Simulation models for a $2^4$ design; the $\beta_i$'s are uniformly distributed and $\varepsilon$ is normally distributed with mean 0 and variance 0.25; $U(a,b)$ denotes a uniform distribution on the interval $(a,b)$; $\varepsilon$ and the $\beta_i$'s are mutually independent* tions with a fixed set of regression coefficients, we randomly picked the coefficients from a uniform distribution in each simulation trial. The *Null* model serves as a baseline where none of the predictor variables has any effect on the mean yield, i.e., the true model is a constant. The *Unif* model has main and interaction effects independently drawn from a uniform distribution on the interval $(-0.25, 0.25)$. The *Hier* model follows the hierarchical ordering principle—its interaction effects are formed from products of main effects that are bounded by 1 in absolute value. Thus higher-order interaction effects are smaller in magnitude than their lower-order parent effects. Finally, the *Exp* model has non-normal errors and variance heterogeneity, with the variance increasing with the mean.

| Name | Simulation model | $\beta$ distribution |
|---|---|---|
| Null | $y = \varepsilon$ | |
| Unif | $y = \beta_1x_1 + \beta_2x_2 + \beta_3x_3 + \beta_4x_4 + \beta_5x_1x_2 + \beta_6x_1x_3 + \beta_7x_1x_4 + \beta_8x_2x_3 + \beta_9x_2x_4 + \beta_{10}x_3x_4 + \beta_{11}x_1x_2x_3 + \beta_{12}x_1x_2x_4 + \beta_{13}x_1x_3x_4 + \beta_{14}x_2x_3x_4 + \beta_{15}x_1x_2x_3x_4 + \varepsilon$ | $U(-1/4, 1/4)$ |
| Exp | $y = \exp(\beta_1x_1 + \beta_2x_2 + \beta_3x_3 + \beta_4x_4 + \varepsilon)$ | $U(-1, 1)$ |
| Hier | $y = \beta_1x_1 + \beta_2x_2 + \beta_3x_3 + \beta_4x_4 + \beta_1\beta_2x_1x_2 + \beta_1\beta_3x_1x_3 + \beta_1\beta_4x_1x_4 + \beta_2\beta_3x_2x_3 + \beta_2\beta_4x_2x_4 + \beta_3\beta_4x_3x_4 + \beta_1\beta_2\beta_3x_1x_2x_3 + \beta_1\beta_2\beta_4x_1x_2x_4 + \beta_1\beta_3\beta_4x_1x_3x_4 + \beta_2\beta_3\beta_4x_2x_3x_4 + \beta_1\beta_2\beta_3\beta_4x_1x_2x_3x_4 + \varepsilon$ | $U(-1, 1)$ |

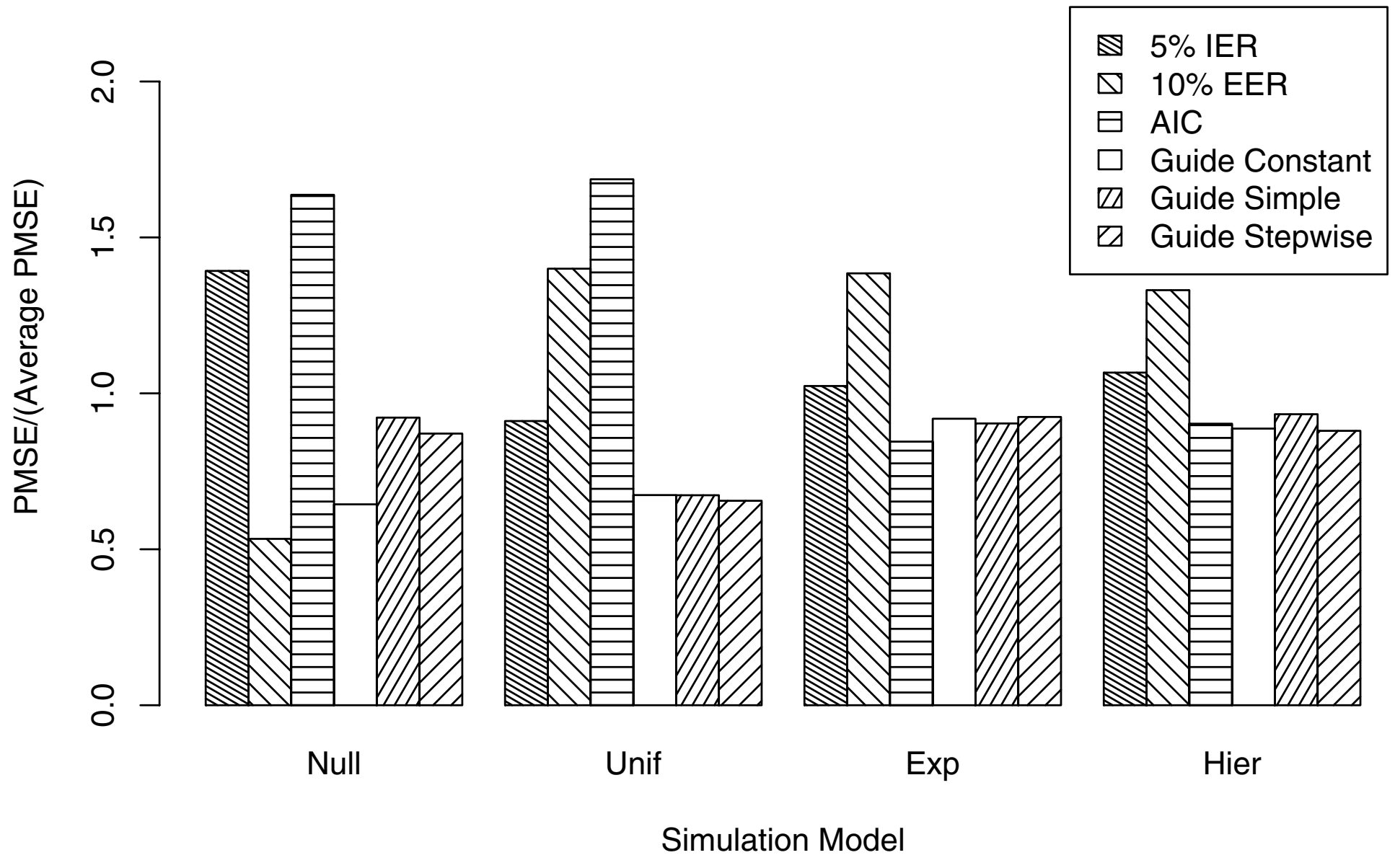

FIG 4. *Barplots of relative PMSE of methods for the four simulation models in Table* 2. *The relative PMSE of a method at a simulation model is defined as its PMSE divided by the average PMSE of the six methods at the same model.*

Ten thousand simulation trials were performed for each model. For each trial, 96 observations were simulated, yielding 6 replicates at each of the 16 factor-level combinations of a $2^4$ design. Each method was applied to find estimates, $\hat{\mu}_i$, of the 16 true means, $\mu_i$, and the sum of squared errors $\sum_1^{16}(\hat{\mu}_i - \mu_i)^2$ was computed. The average over the 10,000 simulation trials gives an estimate of the PMSE of the method. Figure 4 shows barplots of the relative PMSEs, where each PMSE is divided by the average PMSE over the methods. This is done to overcome differences in the scale of the PMSEs among simulation models. Except for a couple of bars of almost identical lengths, the differences in length for all the other bars are statistically significant at the 0.1-level according to Tukey HSD simultaneous confidence intervals.

It is clear from the lengths of the bars for the IER and AIC methods under the *Null* model that they tend to overfit the data. Thus they are more likely than the other methods to identify an effect as significant when it is not. As may be expected, the EER method performs best at controlling the probability of false positives. But it has the highest PMSE values under the non-null situations. In contrast, the three

GUIDE methods provide a good compromise; they have relatively low PMSE values across all four simulation models.

## 4. Unreplicated $2^5$ experiments

If an experiment is unreplicated, we cannot get an unbiased estimate of $\sigma^2$. Consequently, the IER and ERR approaches to model selection cannot be applied. The AIC method is useless too because it always selects the full model. For two-level factorial experiments, practitioners often use a rather subjective technique, due to Daniel [11], that is based on a half-normal quantile plot of the absolute estimated main and interaction effects. If the true effects are all null, the plotted points would lie approximately on a straight line. Daniel's method calls for fitting a line to a subset of points that appear linear near the origin and labeling as outliers those that fall far from the line. The selected model is the one that contains only the effects associated with the outliers.

For example, consider the data from a $2^5$ reactor experiment given in Box, Hunter, and Hunter [1, p. 260]. There are 32 observations on five variables and Figure 5 shows a half-normal plot of the estimated effects. The authors judge that there are only five significant effects, namely, $B, D, E, BD$, and $DE$, yielding the model

$$\hat{y} = 65.5 + 9.75x_2 + 5.375x_4 - 3.125x_5 + 6.625x_2x_4 - 5.5x_4x_5. \tag{4.1}$$

Because Daniel did not specify how to draw the straight line and what constitutes an outlier, his method is difficult to apply objectively and hence cannot be evaluated by simulation. Formal algorithmic methods were proposed by Lenth [16], Loh [17], and Dong [12]. Lenth's method is the simplest. Based on the tables in Wu and Hamada [22, p. 620], the 0.05 IER version of Lenth's method gives the same model

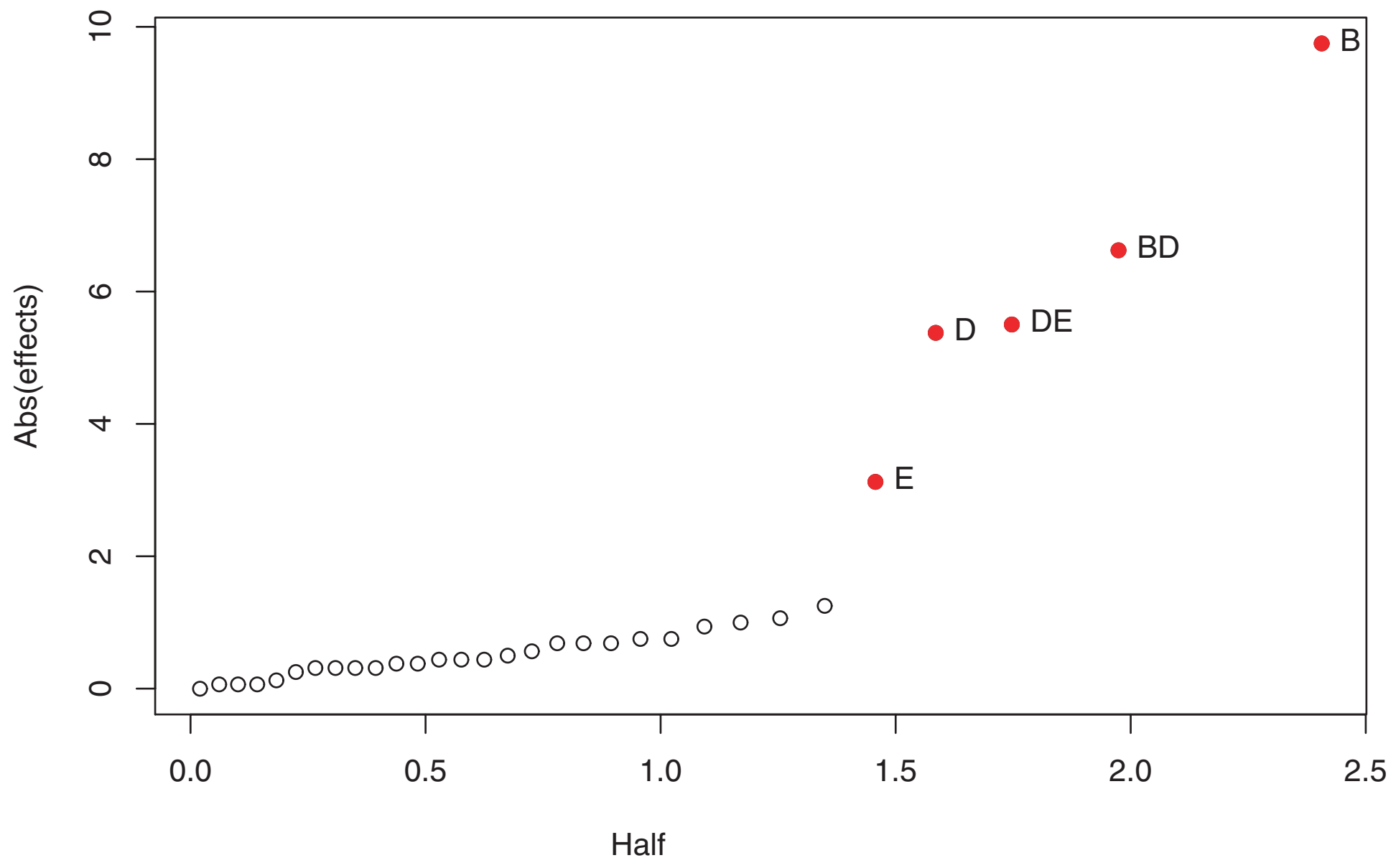

Fig 5. *Half-normal quantile plot of estimated effects from $2^5$ reactor experiment.*

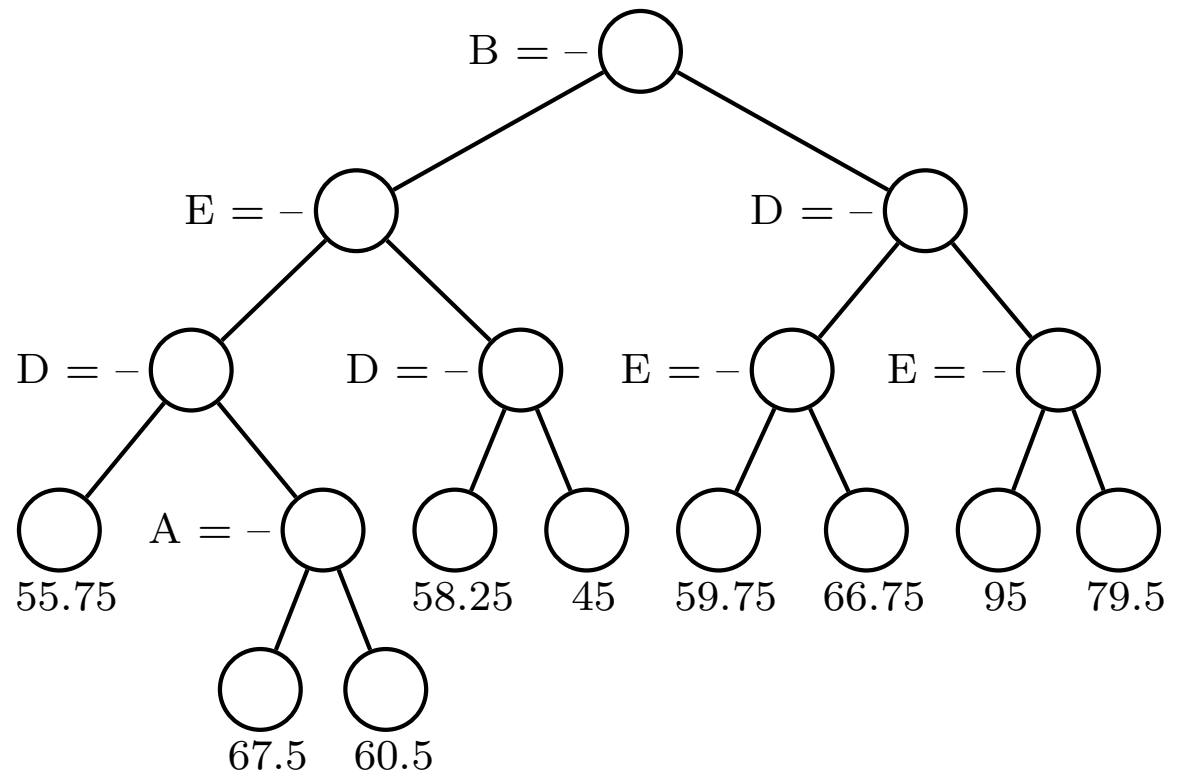

FIG 6. *Piecewise constant GUIDE model for the $2^5$ reactor experiment. The sample y-mean is given beneath each leaf node.*

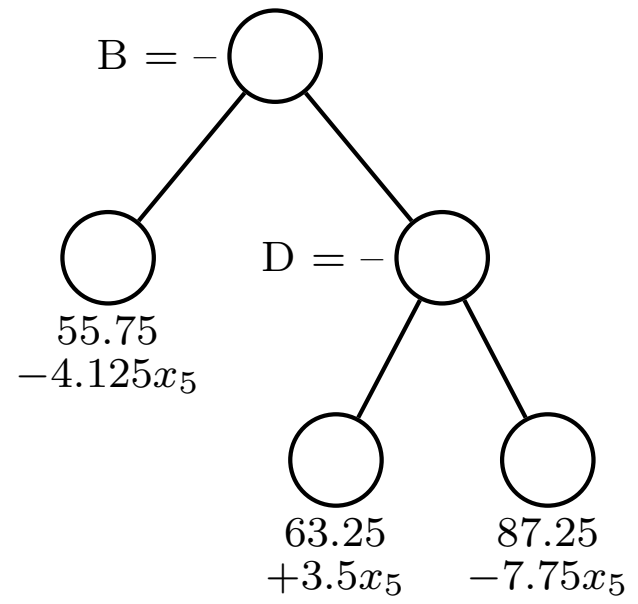

FIG 7. *Piecewise simple linear GUIDE model for the $2^5$ reactor experiment. The fitted equation is given beneath each leaf node.*

as (4.1). The 0.1 EER version drops the $E$ main effect, giving

$$\hat{y} = 65.5 + 9.75x_2 + 5.375x_4 + 6.625x_2x_4 - 5.5x_4x_5. \tag{4.2}$$

The piecewise constant GUIDE model for this dataset is shown in Figure 6. Besides variables $B$, $D$, and $E$, it finds that variable $A$ also has some influence on the yield, albeit in a small region of the design space. The maximum predicted yield of 95 is attained when $B = D = +$ and $E = -$, and the minimum predicted yield of 45 when $B = -$ and $D = E = +$.

If at each node, instead of fitting a constant we fit a best simple linear regression model, we obtain the tree in Figure 7. Factor $E$, which was used to split the nodes at the second and third levels of the piecewise constant tree, is now selected as the best linear predictor in all three leaf nodes. We can try to further simplify the tree structure by fitting a multiple linear regression in each node. The result, shown on the left side of Figure 8, is a tree with only one split, on factor $D$. This model was also found by Cheng and Li [8], who use a method called principal Hessian directions to search for linear functions of the regressor variables; see Filliben and Li [13] for another example of this approach.

We can simplify the model even more by replacing multiple linear regression with stepwise regression at each node. The result is shown by the tree on the right side of Figure 8. It is almost the same as the tree on its left, except that only factors $B$ and $E$ appear as regressors in the leaf nodes. This coincides with the Box, Hunter, and Hunter model (4.1), as seen by expressing the tree model algebraically as

$$\begin{aligned}\hat{y} &= (60.125 + 3.125x_2 + 2.375x_5)(1 - x_4)/2 \\ &\quad + (70.875 + 16.375x_2 - 8.625x_5)(1 + x_4)/2 \\ &= 65.5 + 9.75x_2 + 5.375x_4 - 3.125x_5 + 6.625x_2x_4 - 5.5x_4x_5.\end{aligned} \tag{4.3}$$

An argument can be made that the tree model on the right side of Figure 8 provides a more intuitive explanation of the $BD$ and $DE$ interactions than equation (4.4). For example, the coefficient for the $x_2x_4$ term (i.e., $BD$ interaction) in (4.4) is $6.625 = (16.375 - 3.125)/2$, which is half the difference between the coefficients of the $x_2$ terms (i.e., $B$ main effects) in the two leaf nodes of the tree. Since the root node is split on $D$, this matches the standard definition of the $BD$ interaction as half the difference between the main effects of $B$ conditional on the levels of $D$.

How do the five models compare? Their fitted values are very similar, as Figure 9 shows. Note that every GUIDE model satisfies the heredity principle, because by construction an $n$th-order interaction effect appears only if the tree has $(n + 1)$ levels of splits. Thus if a model contains a cross-product term, it must also contain cross-products of all subsets of those variables.

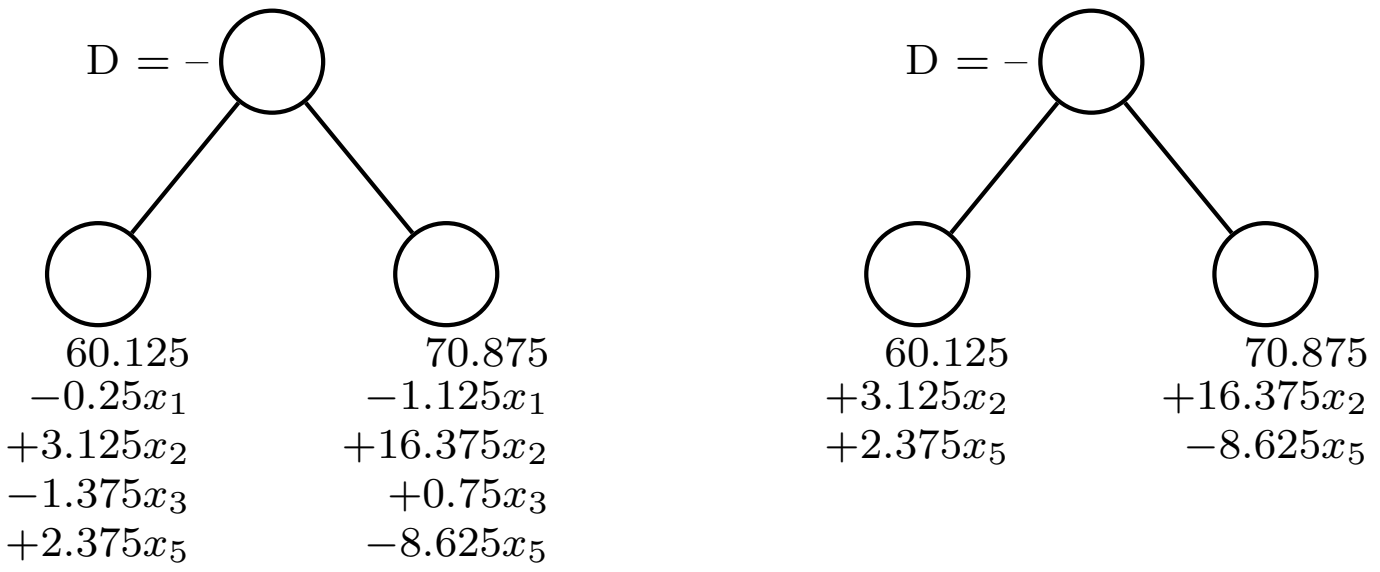

FIG 8. *GUIDE piecewise multiple linear (left) and stepwise linear (right) models.*

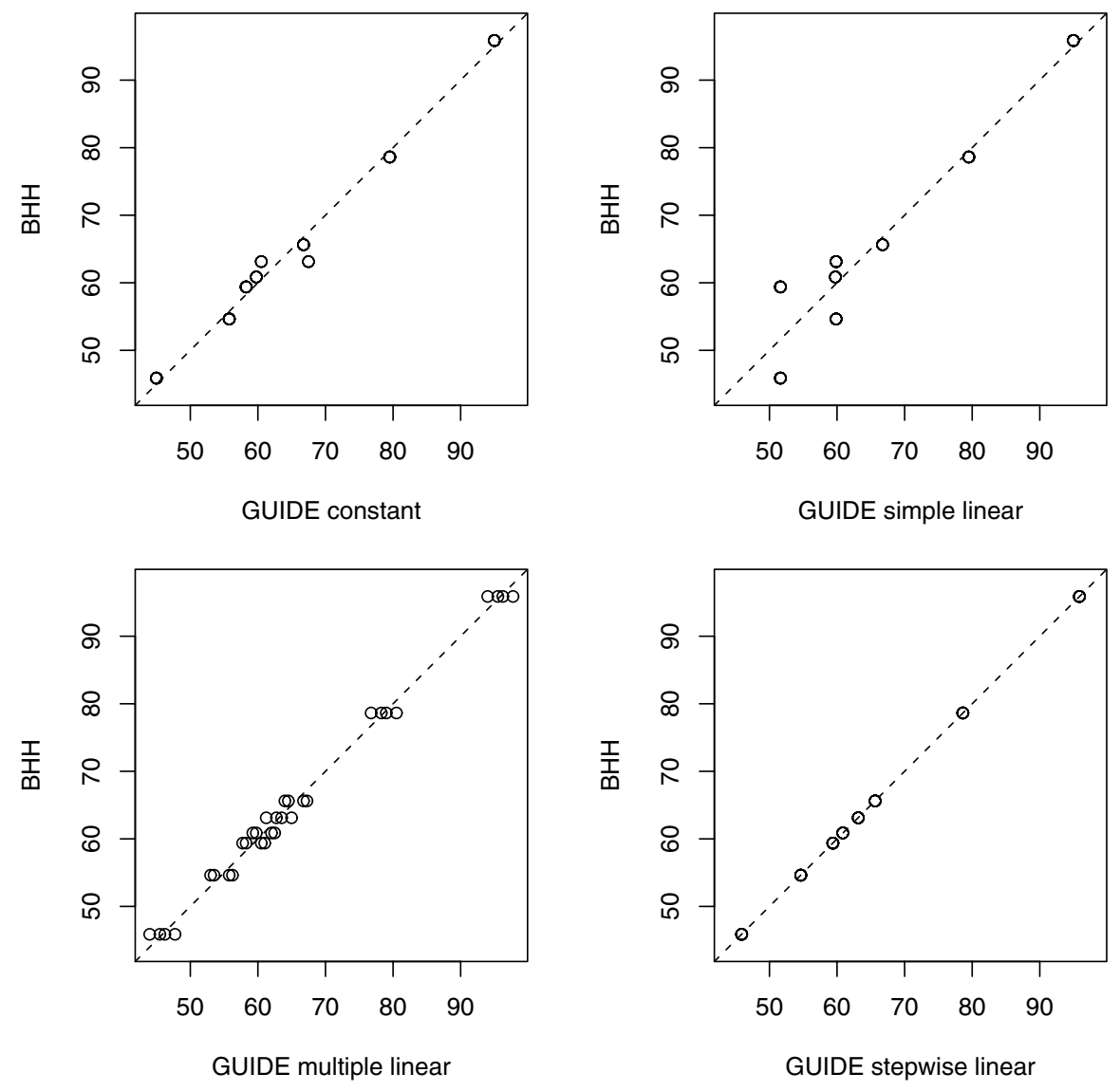

FIG 9. *Plots of fitted values from the Box, Hunter, and Hunter (BHH) model versus fitted values from four GUIDE models for the unreplicated $2^5$ example.*

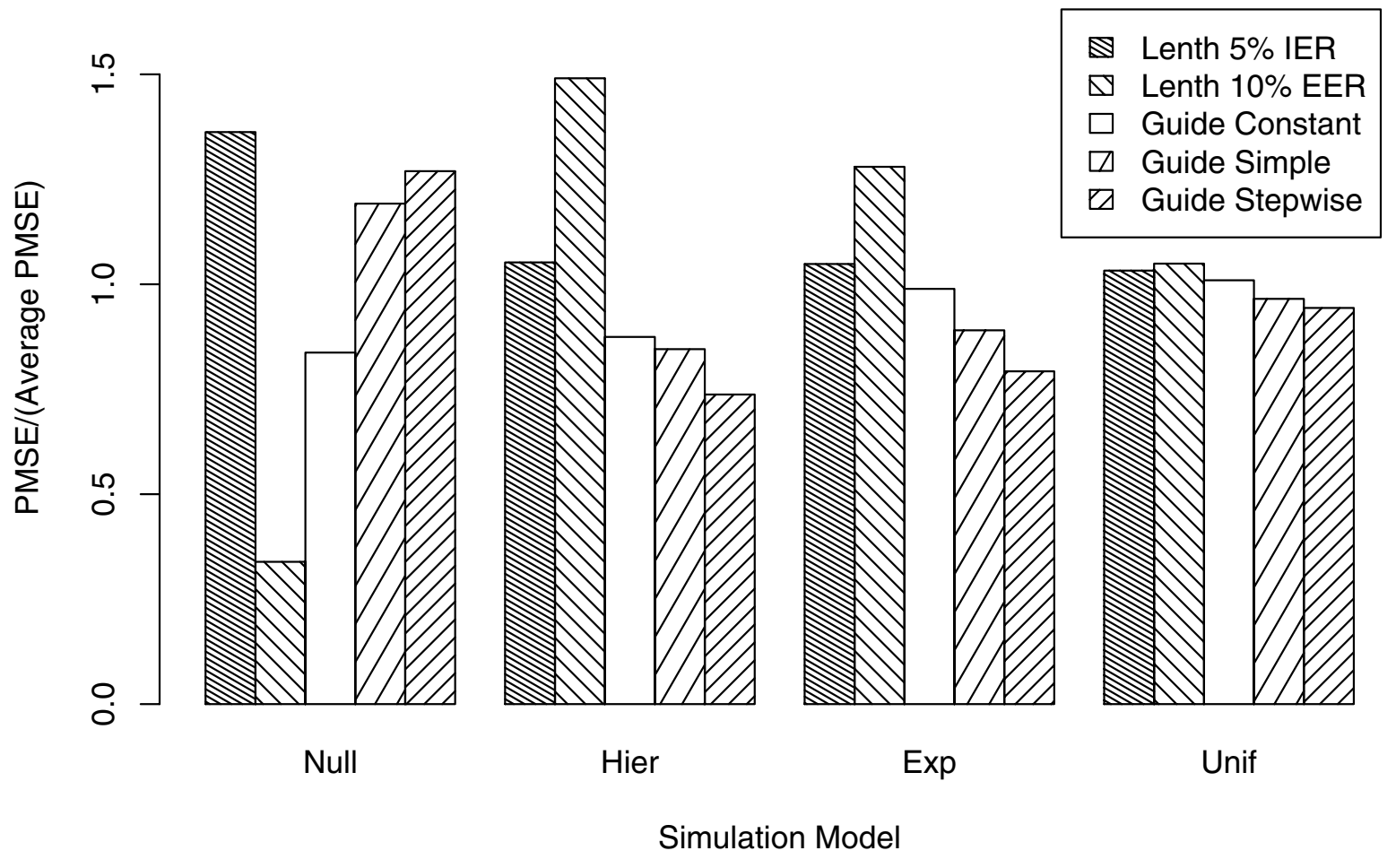

Fig 10. *Barplots of relative PMSEs of Lenth and GUIDE methods for four simulation models. The relative PMSE of a method at a simulation model is defined as its PMSE divided by the average PMSE of the five methods at the same model.*

Figure 10 shows barplots of the simulated relative PMSEs of the five methods for the four simulation models in Table 2. The methods being compared are: (i) Lenth using 0.05 IER, (ii) Lenth using 0.1 EER, (iii) piecewise constant GUIDE, (iv) piecewise best simple linear GUIDE, and (v) piecewise stepwise linear GUIDE. The results are based on 10,000 simulation trials with each trial consisting of 16 observations from an unreplicated $2^4$ factorial. The behavior of the GUIDE models is quite similar to that for replicated experiments in Section 3. Lenth's EER method does an excellent job in controlling the probability of Type I error, but it does so at the cost of under-fitting the non-null models. On the hand, Lenth's IER method tends to over-fit more than any of the GUIDE methods, across all four simulation models.

## 5. Poisson regression

Model interpretation is much harder if some variables have more than two levels. This is due to the main and interaction effects having more than one degree of freedom. We can try to interpret a main effect by decomposing it into orthogonal contrasts to represent linear, quadratic, cubic, etc., effects, and similarly decompose an interaction effect into products of these contrasts. But because the number of products increases quickly with the order of the interaction, it is not easy to interpret several of them simultaneously. Further, if the experiment is unreplicated, model selection is more difficult because significance test-based and AIC-based methods are inapplicable without some assumptions on the order of the correct model.

To appreciate the difficulties, consider an unreplicated $3\times2\times4\times10\times3$ experiment on wave-soldering of electronic components in a printed circuit board reported in Comizzoli, Landwehr, and Sinclair [10]. There are 720 observations and the variables and their levels are:

TABLE 3
*Results from a second-order Poisson loglinear model fitted to solder data*

| Term | Df | Sum of Sq | Mean Sq | F | Pr(>F) |
|---|---|---|---|---|---|
| Opening | 2 | 1587.563 | 793.7813 | 568.65 | 0.00000 |
| Solder | 1 | 515.763 | 515.7627 | 369.48 | 0.00000 |
| Mask | 3 | 1250.526 | 416.8420 | 298.62 | 0.00000 |
| Pad | 9 | 454.624 | 50.5138 | 36.19 | 0.00000 |
| Panel | 2 | 62.918 | 31.4589 | 22.54 | 0.00000 |
| Opening:Solder | 2 | 22.325 | 11.1625 | 8.00 | 0.00037 |
| Opening:Mask | 6 | 66.230 | 11.0383 | 7.91 | 0.00000 |
| Opening:Pad | 18 | 45.769 | 2.5427 | 1.82 | 0.01997 |
| Opening:Panel | 4 | 10.592 | 2.6479 | 1.90 | 0.10940 |
| Solder:Mask | 3 | 50.573 | 16.8578 | 12.08 | 0.00000 |
| Solder:Pad | 9 | 43.646 | 4.8495 | 3.47 | 0.00034 |
| Solder:Panel | 2 | 5.945 | 2.9726 | 2.13 | 0.11978 |
| Mask:Pad | 27 | 59.638 | 2.2088 | 1.58 | 0.03196 |
| Mask:Panel | 6 | 20.758 | 3.4596 | 2.48 | 0.02238 |
| Pad:Panel | 18 | 13.615 | 0.7564 | 0.54 | 0.93814 |
| Residuals | 607 | 847.313 | 1.3959 | | |

1. `Opening:` amount of clearance around a mounting pad (levels 'small', 'medium', or 'large')
2. `Solder:` amount of solder (levels 'thin' and 'thick')
3. `Mask:` type and thickness of the material for the solder mask (levels A1.5, A3, B3, and B6)
4. `Pad:` geometry and size of the mounting pad (levels D4, D6, D7, L4, L6, L7, L8, L9, W4, and W9)
5. `Panel:` panel position on a board (levels 1, 2, and 3)

The response is the number of solder skips, which ranges from 0 to 48.

Since the response variable takes non-negative integer values, it is natural to fit the data with a Poisson log-linear model. But how do we choose the terms in the model? A straightforward approach would start with an ANOVA-type model containing all main effect and interaction terms and then employ significance tests to find out which terms to exclude. We cannot do this here because fitting a full model to the data leaves no residual degrees of freedom for significance testing. Therefore we have to begin with a smaller model and hope that it contains all the necessary terms.

If we fit a second-order model, we obtain the results in Table 3. The three most significant two-factor interactions are between `Opening`, `Solder`, and `Mask`. These variables also have the most significant main effects. Chambers and Hastie [4, p. 10]—see also Hastie and Pregibon [14, p. 217]—determine that a satisfactory model for these data is one containing all main effect terms and these three two-factor interactions. Using set-to-zero constraints (with the first level in alphabetical order set to 0), this model yields the parameter estimates given in Table 4. The model is quite complicated and is not easy to interpret as it has many interaction terms. In particular, it is hard to explain how the interactions affect the mean response.

Figure 11 shows a piecewise constant Poisson regression GUIDE model. Its size is a reflection of the large number of variable interactions in the data. More interesting, however, is the fact that the tree splits first on `Opening`, `Mask`, and `Solder`—the three variables having the most significant two-factor interactions.

As we saw in the previous section, we can simplify the tree structure by fitting a main effects model to each node instead of a constant. This yields the much smaller piecewise main effect GUIDE tree in Figure 12. It has only two splits, first on `Solder` and then, if the latter is thin, on `Opening`. Table 5 gives the regression coefficients in the leaf nodes and Figure 13 graphs them for each level of `Mask` and `Pad` by leaf node.

Table 4
*A Poisson loglinear model containing all main effects and all two-factor interactions involving* `Opening`, `Solder`, *and* `Mask`.

| Regressor | Coef | t | Regressor | Coef | t |
|---|---|---|---|---|---|
| Constant | -2.668 | -9.25 | | | |
| maskA3 | 0.396 | 1.21 | openmedium | 0.921 | 2.95 |
| maskB3 | 2.101 | 7.54 | opensmall | 2.919 | 11.63 |
| maskB6 | 3.010 | 11.36 | soldthin | 2.495 | 11.44 |
| padD6 | -0.369 | -5.17 | maskA3:openmedium | 0.816 | 2.44 |
| padD7 | -0.098 | -1.49 | maskB3:openmedium | -0.447 | -1.44 |
| padL4 | 0.262 | 4.32 | maskB6:openmedium | -0.032 | -0.11 |
| padL6 | -0.668 | -8.53 | maskA3:opensmall | -0.087 | -0.32 |
| padL7 | -0.490 | -6.62 | maskB3:opensmall | -0.266 | -1.12 |
| padL8 | -0.271 | -3.91 | maskB6:opensmall | -0.610 | -2.74 |
| padL9 | -0.636 | -8.20 | maskA3:soldthin | -0.034 | -0.16 |
| padW4 | -0.110 | -1.66 | maskB3:soldthin | -0.805 | -4.42 |
| padW9 | -1.438 | -13.80 | maskB6:soldthin | -0.850 | -4.85 |
| panel2 | 0.334 | 7.93 | openmedium:soldthin | -0.833 | -4.80 |
| panel3 | 0.254 | 5.95 | opensmall:soldthin | -0.762 | -5.13 |

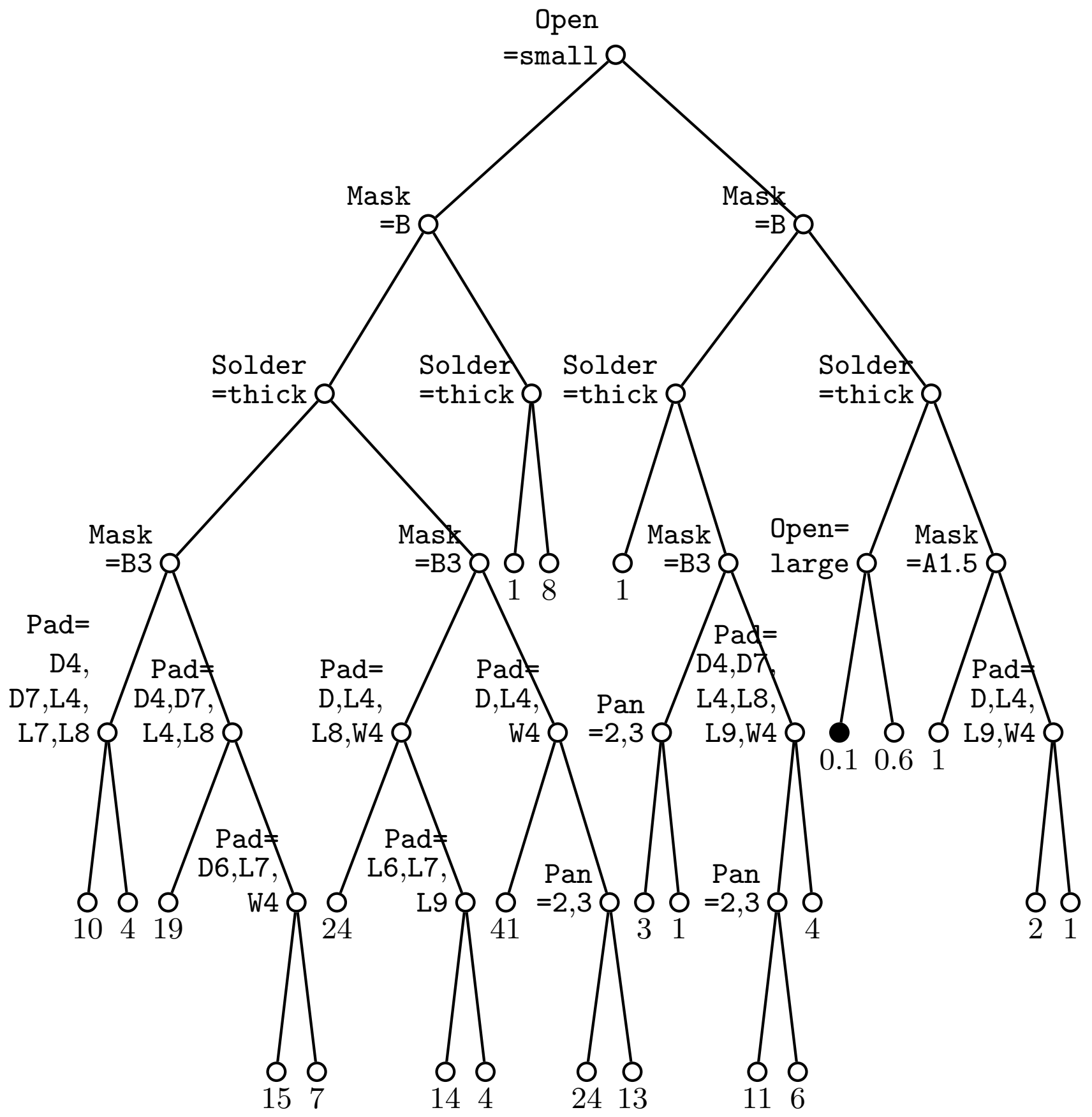

Fig 11. *GUIDE piecewise constant Poisson regression tree for solder data. "Panel" is abbreviated as "Pan". The sample mean yield is given beneath each leaf node. The leaf node with the lowest mean yield is painted black.*

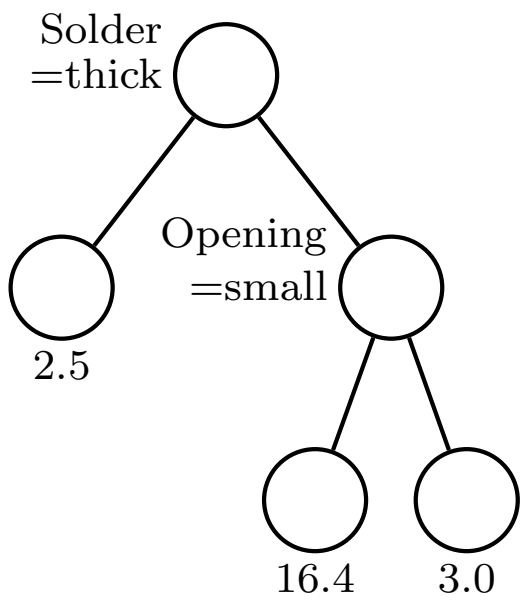

Fig 12. *GUIDE piecewise main effect Poisson regression tree for solder data. The number beneath each leaf node is the sample mean response.*

Table 5
*Regression coefficients in leaf nodes of Figure* 12

| | Solder thick | | Solder thin | | | |
|---|---|---|---|---|---|---|
| | | | Opening small | | Opening not small | |
| Regressor | Coef | t | Coef | t | Coef | t |
| Constant | -2.43 | -10.68 | 2.08 | 21.50 | -0.37 | -1.95 |
| mask=A3 | 0.47 | 2.37 | 0.31 | 3.33 | 0.81 | 4.55 |
| mask=B3 | 1.83 | 11.01 | 1.05 | 12.84 | 1.01 | 5.85 |
| mask=B6 | 2.52 | 15.71 | 1.50 | 19.34 | 2.27 | 14.64 |
| open=medium | 0.86 | 5.57 | aliased | - | 0.10 | 1.38 |
| open=small | 2.46 | 18.18 | aliased | - | aliased | - |
| pad=D6 | -0.32 | -2.03 | -0.25 | -2.79 | -0.80 | -4.65 |
| pad=D7 | 0.12 | 0.85 | -0.15 | -1.67 | -0.19 | -1.35 |
| pad=L4 | 0.70 | 5.53 | 0.08 | 1.00 | 0.21 | 1.60 |
| pad=L6 | -0.40 | -2.46 | -0.72 | -6.85 | -0.82 | -4.74 |
| pad=L7 | 0.04 | 0.29 | -0.65 | -6.32 | -0.76 | -4.48 |
| pad=L8 | 0.15 | 1.05 | -0.43 | -4.45 | -0.36 | -2.41 |
| pad=L9 | -0.59 | -3.43 | -0.64 | -6.26 | -0.67 | -4.05 |
| pad=W4 | -0.05 | -0.37 | -0.09 | -1.00 | -0.23 | -1.57 |
| pad=W9 | -1.32 | -5.89 | -1.38 | -10.28 | -1.75 | -7.03 |
| panel=2 | 0.22 | 2.72 | 0.31 | 5.47 | 0.58 | 5.73 |
| panel=3 | 0.07 | 0.81 | 0.19 | 3.21 | 0.69 | 6.93 |

Because the regression coefficients in Table 5 pertain to conditional main effects only, they are simple to interpret. In particular, all the coefficients except for the constants and the coefficients for `Pad` have positive values. Since negative coefficients are desirable for minimizing the response, the best levels for all variables except `Pad` are thus those not in the table (i.e, whose levels are set to zero). Further, `W9` has the largest negative coefficient among `Pad` levels in every leaf node. Hence, irrespective of `Solder`, the best levels to minimize mean yield are A1.5 `Mask`, large `Opening`, W9 `Pad`, and `Panel` position 1. Finally, since the largest negative constant term occurs when `Solder` is thick, the latter is the best choice for minimizing mean yield. Conversely, it is similarly observed that the worst combination (i.e., one giving the highest predicted mean number of solder skips) is thin `Solder`, small `Opening`, B6 `Mask`, L4 `Pad`, and `Panel` position 2.

Given that the tree has only two levels of splits, it is safe to conclude that four-factor and higher interactions are negligible. On the other hand, the graphs in Figure 13 suggest that there may exist some weak three-factor interactions, such as between `Solder`, `Opening`, and `Pad`. Figure 14, which compares the fits of this model with those of the Chambers-Hastie model, shows that the former fits slightly better.

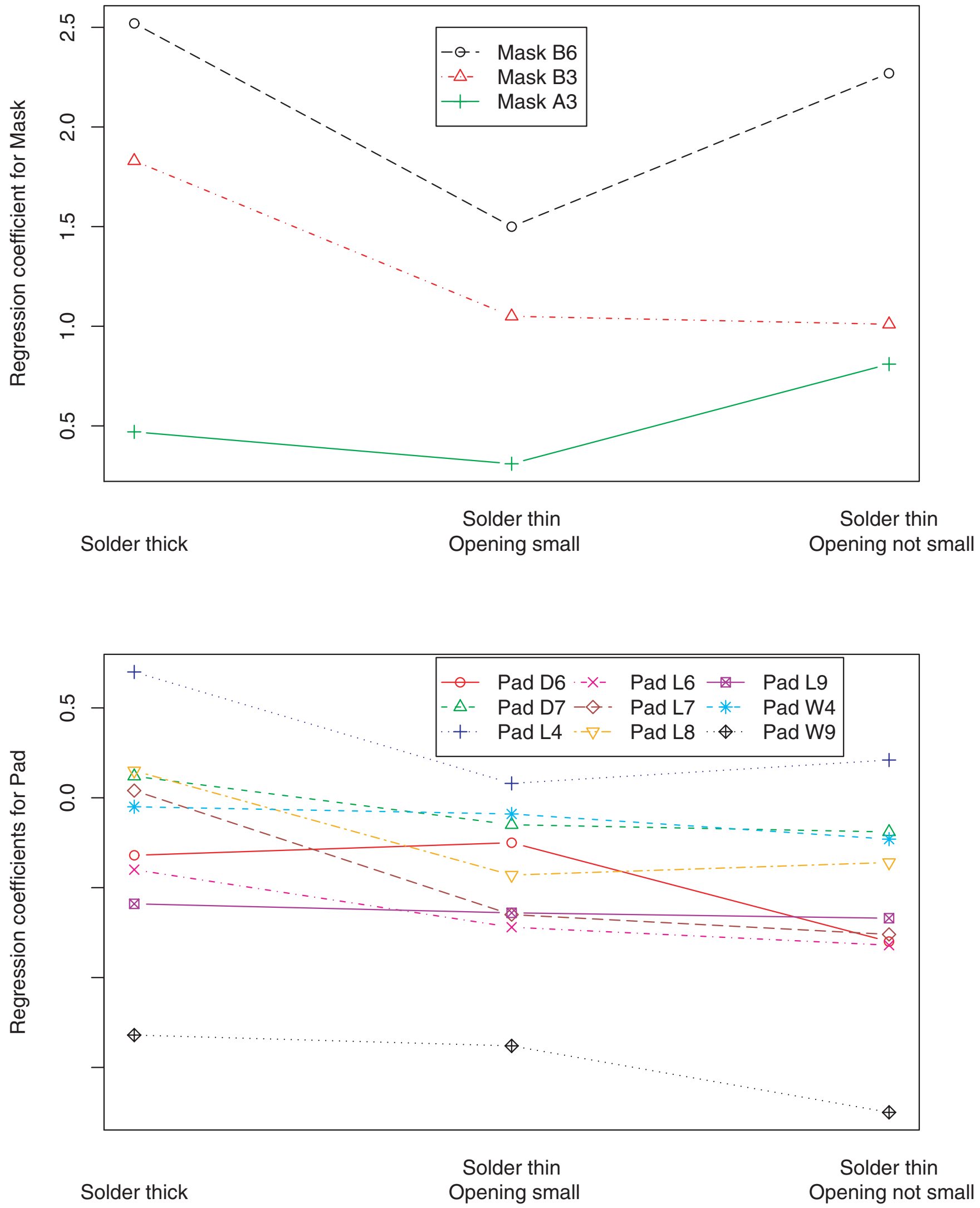

FIG 13. *Plots of regression coefficients for* `Mask` *and* `Pad` *from Table* 5.

## 6. Logistic regression

The same ideas can be applied to fit logistic regression models when the response variable is a sample proportion. For example, Table 6 shows data reported in Collett [9, p. 127] on the number of seeds germinating, out of 100, at two germination temperatures. The seeds had been stored at three moisture levels and three storage temperatures. Thus the experiment is a $2 \times 3 \times 3$ design.

Treating all the factors as nominal, Collett [9, p. 128] finds that a linear logistic regression model with all three main effects and the interaction between moisture level and storage temperature fits the sample proportions reasonably well. The parameter estimates in Table 7 show that only the main effect of storage temperature and its interaction with moisture level are significant at the 0.05 level. Since the storage temperature main effect has two terms and the interaction has four, it takes some effort to fully understand the model.

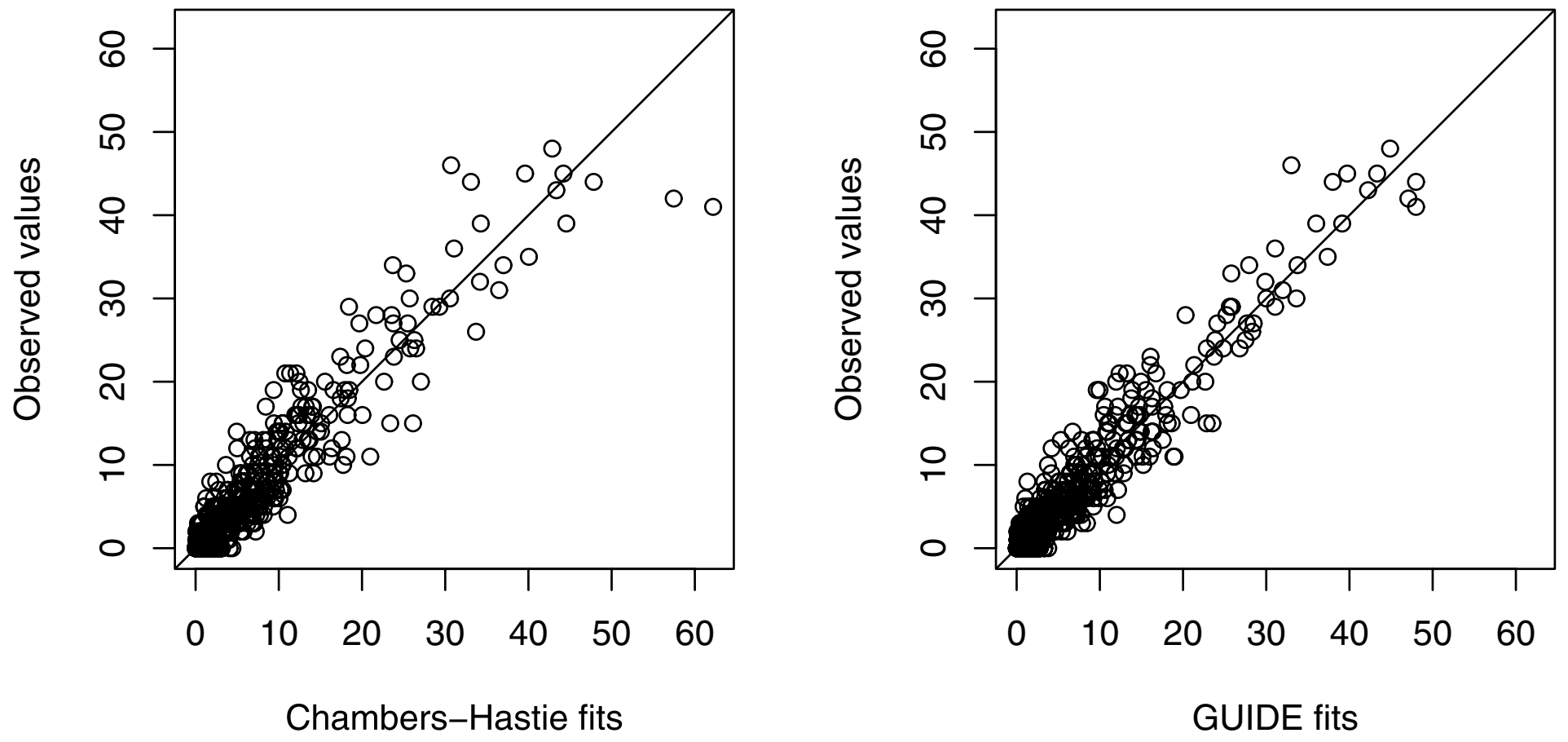

FIG 14. *Plots of observed versus fitted values for the Chambers–Hastie model in Table* 4 *(left) and the GUIDE piecewise main effects model in Table* 5 *(right).*

TABLE 6
*Number of seeds, out of 100, that germinate*

| Germination temp. (°C) | Moisture level | Storage temp. (°C) 21 | 42 | 62 |
|---|---|---|---|---|
| 11 | low | 98 | 96 | 62 |
| 11 | medium | 94 | 79 | 3 |
| 11 | high | 92 | 41 | 1 |
| 21 | low | 94 | 93 | 65 |
| 21 | medium | 94 | 71 | 2 |
| 21 | high | 91 | 30 | 1 |

TABLE 7
*Logistic regression fit to seed germination data using set-to-zero constraints*

| | Coef | SE | $z$ | Pr(> $\|z\|$) |
|---|---|---|---|---|
| (Intercept) | 2.5224 | 0.2670 | 9.447 | < 2e-16 |
| germ21 | -0.2765 | 0.1492 | -1.853 | 0.06385 |
| store42 | -2.9841 | 0.2940 | -10.149 | < 2e-16 |
| store62 | -6.9886 | 0.7549 | -9.258 | < 2e-16 |
| moistlow | 0.8026 | 0.4412 | 1.819 | 0.06890 |
| moistmed | 0.3757 | 0.3913 | 0.960 | 0.33696 |
| store42:moistlow | 2.6496 | 0.5595 | 4.736 | 2.18e-06 |
| store62:moistlow | 4.3581 | 0.8495 | 5.130 | 2.89e-07 |
| store42:moistmed | 1.3276 | 0.4493 | 2.955 | 0.00313 |
| store62:moistmed | 0.5561 | 0.9292 | 0.598 | 0.54954 |

A simple linear logistic regression model, on the other hand, is completely and intuitively explained by its graph. Therefore we will fit a piecewise simple linear logistic model to the data, treating the three-valued storage temperature variable as a continuous linear predictor. We accomplish this with the LOTUS [5] algorithm, which extends the GUIDE algorithm to logistic regression. It yields the logistic regression tree in Figure 15. Since there is only one linear predictor in each node of the tree, the LOTUS model can be visualized through the fitted probability functions shown in Figure 16. Note that although the tree has five leaf nodes, and hence five fitted probability functions, we can display the five functions in three graphs, using solid and dashed lines to differentiate between the two germination temperature levels. Note also that the solid and dashed lines coincide in the middle graph because the fitted probabilities there are independent of germination temperature.

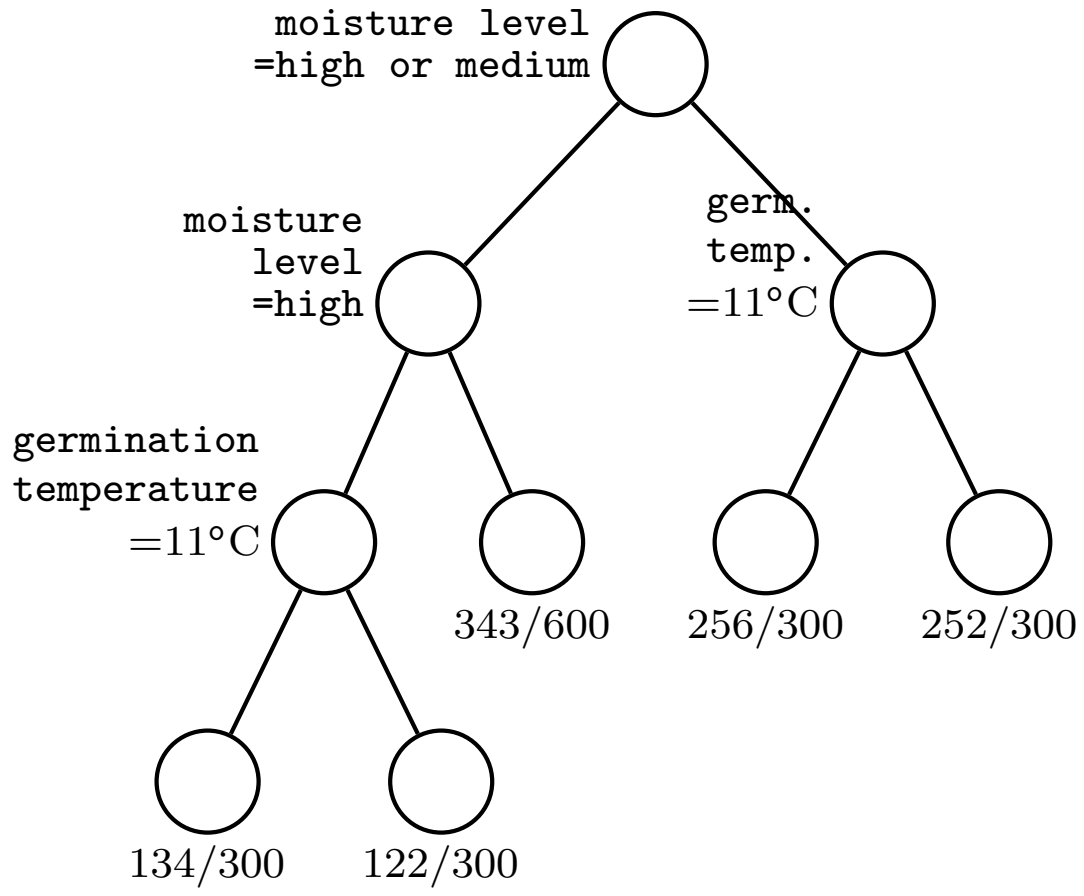

FIG 15. *Piecewise simple linear LOTUS logistic regression tree for seed germination experiment. The fraction beneath each leaf node is the sample proportion of germinated seeds.*

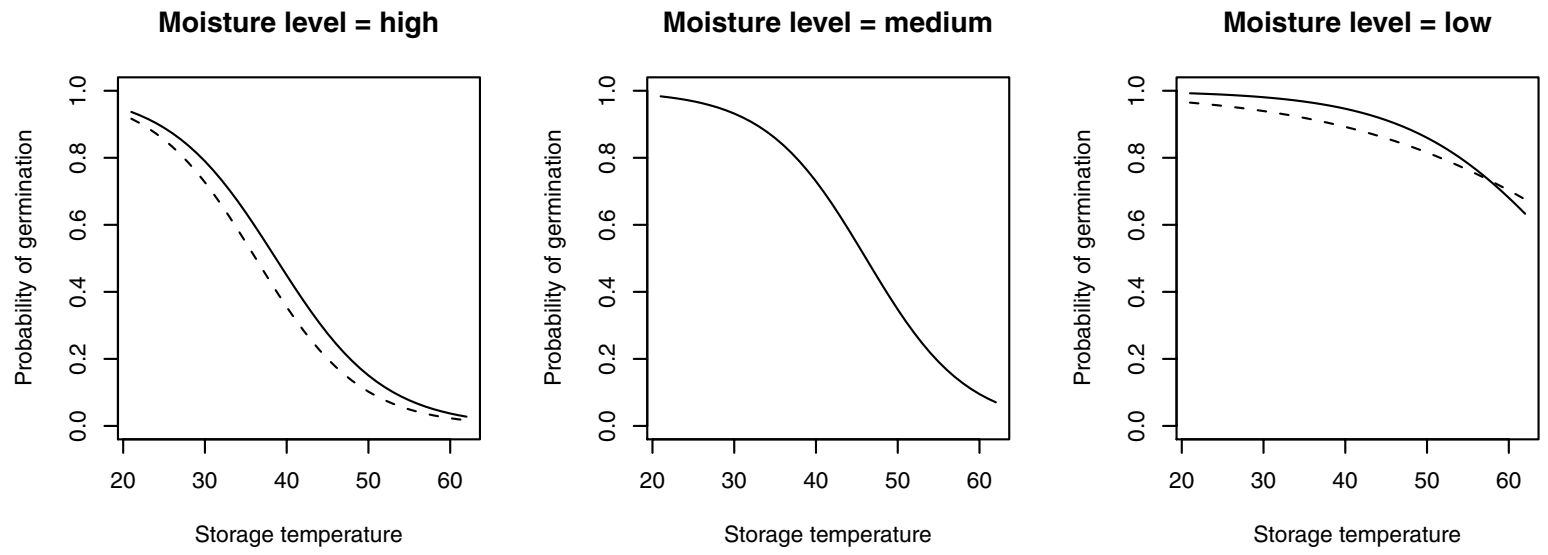

FIG 16. *Fitted probability functions for seed germination data. The solid and dashed lines pertain to fits at germination temperatures of* 11 *and* 21 *degrees, respectively. The two lines coincide in the middle graph.*

The graphs show clearly the large negative effect of storage temperature, especially when moisture level is medium or high. Further, the shapes of the fitted functions for low moisture level are quite different from those for medium and high moisture levels. This explains the strong interaction between storage temperature and moisture level found by Collett [9].

## 7. Conclusion

We have shown by means of examples that a regression tree model can be a useful supplement to a traditional analysis. At a minimum, the former can serve as a check on the latter. If the results agree, the tree offers another way to interpret the main effects and interactions beyond their representations as single degree of freedom contrasts. This is especially important when variables have more than two levels because their interactions cannot be fully represented by low-order contrasts. On the other hand, if the results disagree, the experimenter may be advised to reconsider the assumptions of the traditional analysis. Following are some problems for future study.

1. A tree structure is good for uncovering interactions. If interactions exist, we can expect the tree to have multiple levels of splits. What if there are no interactions? In order for a tree structure to represent main effects, it needs one level of splits for each variable. Hence the complexity of a tree is a sufficient but not necessary condition for the presence of interactions. One way to distinguish between the two situations is to examine the algebraic equation associated with the tree. If there are no interaction effects, the coefficients of the cross-product terms can be expected to be small relative to the main effect terms. A way to formalize this idea would be useful.
2. Instead of using empirical principles to exclude all high-order effects from the start, a tree model can tell us which effects might be important and which unimportant. Here "importance" is in terms of prediction error, which is a more meaningful criterion than statistical significance in many applications. High-order effects that are found this way can be included in a traditional stepwise regression analysis.
3. How well do the tree models estimate the true response surface? The only way to find out is through computer simulation where the true response function is known. We have given some simulation results to demonstrate that the tree models can be competitive in terms of prediction mean squared error, but more results are needed.
4. Data analysis techniques for designed experiments have traditionally focused on normally distributed response variables. If the data are not normally distributed, many methods are either inapplicable or become poor approximations. Wu and Hamada [22, Chap. 13] suggest using generalized linear models for count and ordinal data. The same ideas can be extended to tree models. GUIDE can fit piecewise normal or Poisson regression models and LOTUS can fit piecewise simple or multiple linear logistic models. But what if the response variable takes unordered nominal values? There is very little statistics literature on this topic. Classification tree methods such as CRUISE [15] and QUEST [19] may provide solutions here.
5. Being applicable to balanced as well as unbalanced designs, tree methods can be useful in experiments where it is impossible or impractical to obtain observations from particular combinations of variable levels. For the same reason, they are also useful in response surface experiments where observations are taken sequentially at locations prescribed by the shape of the surface fitted up to that time. Since a tree algorithm fits the data piecewise and hence locally, all the observations can be used for model fitting even if the experimenter is most interested in modeling the surface in a particular region of the design space.